\documentclass{article}
\usepackage{amssymb}
\newtheorem{theorem}{Theorem}[section]

\begin{document}

\title{A new sharp estimate on the dimension of the attractor for the Dirichlet problem of the complex Ginzburg-Landau equation\footnote{ 
E-mail: karan@aegean.gr}}
\date{}
\author{Nikos. I. Karachalios, \\
{\it Department of Mathematics},\\ 
{\it University of the Aegean},\\
{\it Karlovassi, 83200 Samos, Greece}}
\maketitle
\begin{abstract}
Using the improved lower bound on the sum of the eigenvalues of the Dirichlet Laplacian proved by A. D. Melas (Proc. Amer. Math. Soc. \textbf{131} (2003) 631-636), we report a new and sharp estimate for the dimension of the global attractor associated to the complex Ginzburg-Landau equation supplemented with Dirichlet boundary conditions.
\end{abstract}
\section{Introduction}
Let $\Omega\subset\mathbb{R}^n$, $n\geq 1$, be a bounded open set with boundary $\partial\Omega$ and consider the eigenvalues $0<\Lambda_1(\Omega)\leq\Lambda_2(\Omega)\ldots\leq\Lambda_m(\Omega)\leq\ldots$ (repeated with multiplicity) of the Dirichlet Laplacian
\begin{eqnarray}
\label{D1}
-\Delta u&=&\Lambda u,\;\;\mbox{in}\;\;\Omega,\\
u&=&0,\;\;\mbox{on}\;\;\partial\Omega.\nonumber
\end{eqnarray}
A. D. Melas in \cite[Theorem 1, pg. 632 \& pg. 635]{Melas2002}, proved that
\begin{eqnarray}
\label{Toni1}
&&\sum_{j=1}^m\Lambda_{i}(\Omega)\geq\frac{nC_n}{n+2}V(\Omega)^{-\frac{2}{n}}m^{\frac{n+2}{n}}+M_n\frac{V(\Omega)}{I(\Omega)}\;m,\\
&&M_n=\frac{c}{n+2},\;\;\mbox{with $c<(2\pi)^2\omega_n^{-\frac{4}{n}}$, but $c$ independent of $n$}.\nonumber
\end{eqnarray}
Here $\omega_n$ is the volume of the unit ball in $\mathbb{R}^n$, $C_n=(2\pi)^2\omega_n^{-2/n}$ (known as H. Weyl's constant) and $V(\Omega)$ is the volume of $\Omega$. With $I(\Omega)$ we denote the ``moment of inertia'' of $\Omega$, defined as
\begin{eqnarray*}
I(\Omega)=\min_{\alpha\in\mathbb{R}^N}\int_{\Omega}|x-\alpha|^2dx. 
\end{eqnarray*}
Note that if $R:=\left(\frac{V(\Omega)}{\omega_n}\right)^{\frac{1}{n}}$, we have by translating the open set $\Omega$ that 
$$I(\Omega)\geq \int_{B(R)}|x|^2dx=\frac{n\omega_n R^{n+2}}{n+2},$$
see \cite[pg. 635]{Melas2002}.
The estimate (\ref{Toni1}) is an improvement of order $m$ of the estimate of P. Li and S. T. Yau  \cite[Theorem 1, pg. 312]{LiYau83}
\begin{eqnarray}
\label{LiYau1}
\sum_{i=1}^m\Lambda_{i}(\Omega)\geq\frac{n C_n}{n+2}V(\Omega)^{-\frac{2}{n}}m^{\frac{n+2}{n}}.
\end{eqnarray}
Consider  next the Dirichlet initial-boundary value problem for the complex Ginzburg-Landau equation 
\begin{eqnarray}
\label{GL}
\partial_t u-(\lambda+\mathrm{i}\alpha)\Delta u &+&(\kappa+\mathrm{i}\beta)|u|^2u-\gamma u=0,\;\;\mbox{in}\;\;\Omega,t>0,\\
\label{GL2}
u(x,0)&=&u_0(x),\;\;\mbox{for}\;\;x\in\Omega,\\
\label{GL2'}
u(x,t)&=&0\;\;\mbox{in}\;\;\partial\Omega,\;t>0,
\end{eqnarray}
The parameters $\lambda,\alpha,\kappa,\beta,\gamma$ are real and $\lambda,\kappa>0$ (rendering (\ref{GL}) dissipative.  It is well known (see J. M. Ghidaglia \& M. Heron \cite{GhidHer87}, R. Temam \cite{RTem88}) that (\ref{GL})-(\ref{GL2'}) defines a semiflow $\mathcal{S}(t):L^2(\Omega)\rightarrow L^2(\Omega)$ possessing a global attractor $\mathcal{A}$ of finite Hausdorff (and fractal) dimension. We are interested in the case of non-trivial dynamics and we assume that 
\begin{eqnarray}
\label{nontri}
\lambda\Lambda_1<\gamma.
\end{eqnarray}
The case 
\begin{eqnarray}
\label{tri}
\gamma\leq\lambda\Lambda_1,
\end{eqnarray}
is leading to trivial dynamics in the sense $||u(t)||_{L^2(\Omega)}^2\rightarrow 0$ as $t\rightarrow\infty$ and $\mathcal{A}=\{0\}$ (see \cite{RTem88}).

By using the improved lower bound (\ref{Toni1}), we report in this note a new and sharp estimate (see Theorem \ref{main}) on the Hausdorff dimension $\mathrm{dim}_H\mathcal{A}$. The estimate reveals that the Hausdorff (and fractal) dimension of the global attractor associated to (\ref{GL})-(\ref{GL2'}) can be in fact considerably smaller than the one predicted by the known estimates (see also \cite{GhidHer87, RTem88}). We mention that (\ref{Toni1}) can possibly produce interesting improvements and observations on the Hausdorff dimension of global attractors associated to the Dirichlet problem  for other dissipative pde's,  e.g. reaction diffusion equations and damped and driven semilinear wave equations  or Schr\"odinger equations.  An example has been reported in \cite{NKNZ08} concerning Allen-Cahn type parabolic equations. 
\section{ A sharp estimate on the Hausdorff dimension of $\mathcal{A}$.}
For simplicity reasons we denote by $H^1_0(\Omega)$ and $L^2(\Omega)$ the complexified spaces  $\mathbb{H}^1_0(\Omega)\equiv H^1_0(\Omega)^2$ and $\mathbb{L}^2(\Omega)\equiv L^2(\Omega)^2$ for the complex field $u=u_1+\mathrm{i}u_2$. The operator $-(\lambda+\mathrm{i}\alpha)\Delta :H^2(\Omega)\cap H^1_0(\Omega)\rightarrow L^2(\Omega)$ can be rewritten as
\begin{equation}
\label{matrixA}
-\mathbf{M}\Delta,\;\;\mathbf{M}=
\left(
\begin{array}{cc}
\lambda&-\alpha\\
\alpha&\lambda
\end{array}
\right).
\end{equation}
Furthermore, $(\mathbf{M}+\mathbf{M}^*)/2=\lambda\mathbf{I}$. Denoting by $\lambda_j(\Omega)$ the eigenvalues of the Dirichlet eigenvalue problem for the operator $-\mathbf{I}\Delta$ in $\Omega$, it holds that $\left\{\lambda_j(\Omega)\right\}_{j=1}^{\infty}=\left\{\Lambda_{j}(\Omega),\Lambda_{j}(\Omega)\right\}_{j=1}^{\infty}$. By using (\ref{Toni1}) and the reduction-convexity argument of V. V. Chepyzhov and A. A. Ilyin \cite[Theorem 3.1, pg. 816]{Chep99}  we get that
\begin{eqnarray}
\label{c01}
\sum_{j=1}^m\lambda_j\geq \frac{2^{-\frac{2}{n}}n C_n}{(n+2)}V(\Omega)^{-\frac{2}{n}}m^{\frac{n+2}{n}}+M_n\frac{V(\Omega)}{I(\Omega)}m.
\end{eqnarray}
The first variation equation of (\ref{GL}), $$\partial_tU-F'(u)U=0,$$ is equivalent to
\begin{eqnarray}
\label{c02}
\partial_tU-(\lambda+\mathrm{i}\alpha)\Delta U+(\kappa+\mathrm{i}\beta)\left\{|u|^2U+2u\mathrm{Re}(\overline{u}U\right\}-\gamma U=0,
\end{eqnarray}
and it is supplemented with the boundary and initial conditions
\begin{eqnarray}
\label{c03}
U=0\;\;\mbox{on}\;\;\partial\Omega\;\;\mbox{and}\;\;U(0)=\xi\in L^2(\Omega).
\end{eqnarray} 
Problem (\ref{c02})-(\ref{c03}) has a unique solution $U\in L^2([0, T];H^1_0(\Omega))\cap L^{\infty}([0,T];L^2(\Omega))$, for all $T>0$. The function $u_0\rightarrow\mathcal{S}(t)u_0$ is Fr\'{e}chet differentiable and the differential is $L(t,u_0):\xi\in L^2(\Omega)\rightarrow U(t)\in L^2(\Omega)$. These fundamental results are from \cite{GhidHer87} and \cite{RTem88}. The sum $\mu_1+\ldots+\mu_m$ of the first $m$ global Lyapunov exponents is bounded by the quantity
\begin{eqnarray*}
q_m=\limsup_{t\rightarrow\infty}\sup_{||\xi_i||_{L^2(\Omega)}\leq 1}\frac{1}{t}
\int_0^t\mathrm{Re}\mathrm{Tr}[F'(S(\tau)u_0)\circ \mathcal{Q}_m(\tau)]d\tau,
\end{eqnarray*}
where $\mathcal{Q}_m(t)$ is the orthogonal projection in $L^2(\Omega)$ onto $\mathrm{span}\left\{U_1(t),\ldots,U_m(t)\right\}$. The vectors $U_{i}(t)$, $i=1,\ldots,m$ are $m$-solutions of (\ref{c02})-(\ref{c03}) starting from the initial conditions $U_i(0)=\xi_i\in L^2(\Omega)$ and the $m$-dimensional volume $|U_1(t)\wedge U_2(t)\ldots\wedge U_m(t)|$ of the infinitesimal parallelepiped spanned by $U_i(t)$ is given by
\begin{eqnarray}
\label{c04}
|U_1(t)\wedge\Phi_2(t)\ldots\wedge U_m(t)|&=&|\xi_1\wedge\ldots\wedge\xi_m|\exp\int_{0}^{t}\mathrm{Re}\mathrm{Tr}[F'(\mathcal{S}(\tau)u_0)\circ\mathcal{Q}_m(\tau)]ds.
\end{eqnarray}
Finding $m$ such that $q_m<0$, Constantin-Foias-Temam theory \cite{CFT1985,RTem88} implies the exponential decay of the $m$-volume element. We fix $\tau$ for the time being and we consider an orthonormal basis $\{\phi_j(\tau)\}_{j=1}^{\infty}$ of $L^2(\Omega)$ with $$\mathrm{span}\left\{\phi_1(\tau),\ldots,\phi_m(\tau)\right\}=\mathrm{span}\left\{U_1(\tau),\ldots,U_m(\tau)\right\}=\mathrm{span}\mathcal{Q}_m(\tau)L^2(\Omega).$$
Since $U_j(\tau)\in H_0^1(\Omega)$, for all $j\in\mathbb{N}$ and almost all $\tau>0$, we have $\phi_j(\tau)\in H_0^1(\Omega)$ for all $j\in\mathbb{N}$ and almost all $\tau>0$. Then $\mathcal{Q}_m(\tau)\phi_{j}(\tau)=\phi_{j}(\tau)$ if $j\leq m$ and $\mathcal{Q}_m(\tau)\phi_{j}=0$ otherwise. 

We will follow in part the proof of \cite[Theorem 7.1, pg. 458]{RTem88} which will be modified as follows: Starting from the equation 
\begin{eqnarray*}
\mathrm{Re}\mathrm{Tr}[F'(\mathcal{S}(\tau)u_0)\circ\mathcal{Q}_m(\tau)=\sum_{j=1}^{\infty}\mathrm{Re}\left(F'(u(\tau))\circ Q_m(\tau)\phi_j(\tau),\phi_j(\tau)\right)_{L^2(\Omega)}=\sum_{j=1}^{m}\mathrm{Re}\left(F'(u(\tau))\phi_j(\tau),\phi_j(\tau)\right)_{L^2(\Omega)},
\end{eqnarray*}
it is found (cf. \cite[eq. (7.9), pg. 758]{RTem88}) that
\begin{eqnarray}
\label{c1}
\sum_{j=1}^m(F'(u)\phi_j,\phi_j)_{L^2(\Omega)}=-\lambda\sum_{j=1}^m||\phi_j||^2_{H_0^1(\Omega)}+2|\beta|\int_{\Omega}|u|^2\rho dx+\gamma m,
\end{eqnarray}
where
\begin{eqnarray}
\label{c2}
\rho=\rho(x,\tau)=\sum_{j=1}^m|\phi_j(x,\tau)|^2,\;\;\mbox{for almost all}\;x,\tau.
\end{eqnarray}
In the case of the Dirichlet boundary conditions we have the Lieb-Thirring inequality
\begin{eqnarray}
\label{c3}
\int_{\Omega}\rho^{\frac{n+2}{n}}dx\leq C_*\sum_{j=1}^m||\phi_j||^2_{H_0^1(\Omega)}.
\end{eqnarray}
Using (\ref{c3}) and H\"older and Young's inequalities, the estimate (\ref{c1}) becomes
\begin{eqnarray}
\label{c5}
\sum_{j=1}^m(F'(u)\phi_j,\phi_j)_{L^2(\Omega)}&\leq&  -\lambda\sum_{j=1}^m||\phi_j||^2_{H_0^1(\Omega)}
+2|\beta|\,||u||^2_{L^{n+2}(\Omega)}||\rho||_{L^{\frac{n+2}{n}}(\Omega)}+\gamma m\nonumber\\
&\leq&
-\lambda\sum_{j=1}^m||\phi_j||^2_{H_0^1(\Omega)}+2c_1|\beta|\,||u||^2_{L^{n+2}(\Omega)}\left\{\sum_{j=1}^m||\phi_j||^2_{H_0^1(\Omega)}\right\}^{\frac{n}{n+2}}+\gamma m,\;\;c_1=C_*^{\frac{n}{n+2}},\nonumber\\
&\leq&
-\frac{\lambda}{2}\sum_{j=1}^m||\phi_j||^2_{H_0^1(\Omega)}
+c_2|\beta|^{\frac{n+2}{2}}\lambda^{-\frac{n}{2}}||u||_{L^{n+2}(\Omega)}^{n+2}+\gamma m,\\
&&\mbox{with}\;\;c_2=2\left(\frac{2}{n+2}\right)^{\frac{n+2}{2}}(nC_*)^{\frac{n}{2}},\nonumber
\end{eqnarray}
aiming to use (\ref{c01}) instead of the property $||\rho||_{L^1(\Omega)}=m\leq V(\Omega)^{\frac{2}{n+2}}||\rho||_{L^{\frac{n+2}{n}}(\Omega)}$. Indeed, since by \cite[Lemma VI.2.1,pg. 390]{RTem88}
$$\sum_{j=1}^m||\phi_j||_{H^1_0(\Omega)}^2\geq \sum_{j=1}^m\lambda_j,$$
by inserting in (\ref{c5}) the estimate (\ref{c01}), we deduce that 
\begin{eqnarray}
\label{c6}
\sum_{j=1}^m(F'(u)\phi_j,\phi_j)_{L^2(\Omega)}\leq-\frac{\lambda n C_n}{2[2V(\Omega)]^{\frac{2}{n}}(n+2)}m^{\frac{n+2}{n}}+\left[\gamma-\frac{\lambda}{2}\frac{M_n V(\Omega)}{I(\Omega)}\right]m
+c_2|\beta|^{\frac{n+2}{2}}\lambda^{-\frac{n}{2}}||u(\tau)||_{L^{n+2}(\Omega)}^{n+2}.
\end{eqnarray}
Let us note first, that the case   
\begin{eqnarray*}
0<\gamma\leq \frac{\lambda}{2}\frac{M_nV(\Omega)}{I(\Omega)},
\end{eqnarray*}
corresponds to the case of trivial dynamics: from (\ref{Toni1}) 
\begin{eqnarray*}
\gamma\leq\frac{\lambda}{2}\frac{M_nV(\Omega)}{I(\Omega)}<\lambda\frac{nC_n}{n+2}V(\Omega)^{-\frac{2}{n}}+\lambda \frac{M_nV(\Omega)}{I(\Omega)}\leq\lambda\Lambda_1.
\end{eqnarray*}
which is condition (\ref{tri}). Thus, under assumption (\ref{nontri}) it is obvious that  
\begin{eqnarray}
\label{c8}
\frac{\lambda}{2}\frac{M_nV(\Omega)}{I(\Omega)}<\gamma. 
\end{eqnarray}
With (\ref{c8}) in hand, we proceed from (\ref{c6}) and Young's inequality to
\begin{eqnarray}
\label{c9}
\sum_{j=1}^m(F'(u)\phi_j,\phi_j)_{L^2(\Omega)}\leq  -\frac{1}{4}\frac{\lambda n C_n}{[2V(\Omega)]^{\frac{2}{n}}(n+2)}m^{\frac{n+2}{n}}&+&\frac{2^{n+2}V(\Omega)}{(n+2)(\lambda C_n)^{\frac{n}{2}}}\left[\gamma-\frac{\lambda}{2}\frac{M_n V(\Omega)}{I(\Omega)}\right]^{\frac{n+2}{2}}\nonumber\\
&+&c_2|\beta|^{\frac{n+2}{2}}\lambda^{-\frac{n}{2}}
||u(\tau)||_{L^{n+2}(\Omega)}^{n+2}.
\end{eqnarray}
Hence, the number $q_m$ is estimated as
\begin{eqnarray}
\label{c10}
q_m\leq -Am^{\frac{n+2}{n}}+B,\;\;A=\frac{1}{4}\frac{\lambda n C_n}{[2V(\Omega)]^{\frac{2}{n}}(n+2)},
\end{eqnarray}
where now the constant $B$ is given by
\begin{eqnarray}
\label{c4b}
B&=&\frac{2^{n+2}V(\Omega)}{(n+2)(\lambda C_n)^{\frac{n}{2}}}\left[\gamma-\frac{\lambda}{2}\frac{M_n V(\Omega)}{I(\Omega)}\right]^{\frac{n+2}{2}}+
c_2|\beta|^{\frac{n+2}{2}}\lambda^{-\frac{n}{2}}\delta,\\
\delta&:=&\limsup_{t\rightarrow\infty}\sup_{u_0\in\mathcal{A}}\frac{1}{t}\int_{0}^{t}||S(\tau)u_0||_{L^{n+2}(\Omega)}^{n+2}d\tau.
\end{eqnarray}
Explicit estimates for $\delta$ are given in \cite[Remark 7.1, pg. 460]{RTem88}. It is known from \cite[Corollary 2.2, pg. 815]{Chep99} that if $q_m\leq f(m)$, $m=1,2,\ldots,$ where $f(x)$ is a concave function of the continuous variable $x$, and $p(x^*)=0$, then $\mathrm{dim}_H\mathcal{A}\leq x^*$. In our case  
$f(x)=-Ax^{\frac{n+2}{n}}+B$ for $x>0$. We may summarize our observations in
\begin{theorem}
\label{main}
Consider the global atrractor $\mathcal{A}$ of the semiflow $\mathcal{S}(t):L^2(\Omega)\rightarrow L^2(\Omega)$ defined by the Ginzburg-Landau equation supplemented with the Dirichlet boundary conditions with $\lambda>0$, $\kappa>0$ and the parameter $\gamma>0$ satisfying  (\ref{nontri}). Then
$$\mathrm{dim}_H\mathcal{A}\leq \left(\frac{B}{A}\right)^{\frac{n}{n+2}}:=d^*.$$
The constants $A,B>0$, are given by (\ref{c10}) and (\ref{c4b}) respectively. 
\end{theorem}
\vspace{0.1cm}
We remark that by the improved results of \cite{Chep99}, the fractal dimension of $\mathcal{A}$ satisfies $\mathrm{dim}_F\mathcal{A}\leq d^*$. Theorem \ref{main} actually indicates (due to the appearance of the term $-\lambda M_nV(\Omega)/2I(\Omega)$ in (\ref{c4b})) that the number of degrees of freedom for the Ginzburg-Landau semiflow (\ref{GL})-(\ref{GL2'}) when $\lambda\Lambda_1<\gamma$, can be significantly lower than the one predicted by the known upper bounds (e.g. making use of an estimate of the form (\ref{LiYau1})). See also \cite{bab90,GhidHer87,RTem88}.  
\bibliographystyle{amsplain}

\begin{thebibliography}{10}

\bibitem{bab90} {A. V. Babin and M. I. Vishik}, Attractors of evolution equations. Studies in Mathematics and its Applications, 25. North-Holland Publishing Co., Amsterdam, 1992.
%
\bibitem{Chep99} V. V. Chepyzov, A. A. Ilyin, {\em A note on fractal dimension of attractors of dissipative dynamical systems}, Nonlinear Anal. \textbf{44} (6) (2001), 811--819.
%
\bibitem{CFT1985} P. Constantin, C. Foias and R. Temam \emph{Attractors representing turbulent flows}, Memoirs of the AMS, Vol. 53, no. 314, (1985)
%
\bibitem{GhidHer87} J. M. Ghidaglia and B. H\'{e}ron, \emph{Dimension of the attractor associated to the Ginzburg-Landau equation}, Physica D \textbf{28} (1987), 282-304.
%
\bibitem{NKNZ08} N. I. Karachalios and N. B. Zographopoulos \emph{A sharp estimate and change on the dimension of the attractor for Allen-Cahn equations}. Preprint, arXiv:0802.1860v2 [math. AP].
%
\bibitem{LiYau83} P. Li and S. T. Yau, \emph{On the Schr\"odinger equation and the eigenvalue problem}. Comm. Math. Phys. \textbf{88} (1983), no. 3, 309--318.
%
\bibitem{Melas2002} A. Melas, \emph{A lower bound of sums of eigenvalues of the Laplacian}, Proc. Amer. Math. Society \textbf{131} (2) 631-636.

\bibitem{RTem88} R. Temam, {\em Infinite-Dimensional Dynamical
Systems in Mechanics and Physics, 2nd edition}, Springer-Verlag,
New York, 1997.
\end{thebibliography}

\end{document}